\documentclass[11pt,a4paper]{article}

\usepackage{amsmath,amssymb}

\usepackage{color}

\newtheorem{dfn}{Definition}[section]
\newtheorem{thm}[dfn]{Theorem}
\newtheorem{rmk}[dfn]{Remark}
\newtheorem{lem}[dfn]{Lemma}

\newcommand{\midd}{\,|\,}
\newcommand{\W}{\mathcal{W}}
\newcommand{\p}{\partial}
\newcommand{\e}{\varepsilon}
\newcommand{\nn}{\nonumber}
\newcommand{\C}{\mathbb{C}}

\newcommand{\Z}{\mathbb{Z}}

\newcommand{\g}{\mathfrak{g}}

\newcommand{\h}{\mathfrak{h}}
\newcommand{\hh}{\hat{\h}}
\newcommand{\F}{\mathcal{F}}

\newenvironment{prf}{\noindent \textit{Proof} \ }{\hfill $\Box$}
\newenvironment{prfn}[1]{\vskip 1em \noindent \textit{Proof of #1} \ }{\hfill $\Box$}

\begin{document}
\title{Uniqueness Theorem of $\W$-Constraints for Simple Singularities}
\author{Si-Qi Liu, Di Yang, Youjin Zhang\\
\small Department of Mathematical Sciences, Tsinghua University, \\
\small Beijing, 100084, China}
\maketitle

\begin{abstract}
In a recent paper \cite{BM}, Bakalov and Milanov proved that the total descendant potential of a simple singularity satisfies the $\W$-constraints, which come from the $\W$-algebra 
of the lattice vertex algebra associated to the root lattice of this singularity and a twisted module of
this vertex algebra. 
In the present paper, we prove that the solution of these $\W$-constraints is unique up
to a constant factor, as conjectured by Bakalov and Milanov in their paper.
\end{abstract}

\section{Introduction}

In 1991, Witten proposed his celebrated conjecture on the relation between the partition function of
the two-dimensional topological gravity and the Korteweg-de Vries (KdV) integrable hierarchy \cite{W1}.
This conjecture has two equivalent versions. Let us denote the partition function of the two-dimensional topological gravity by $\tau$,
then the first version of the Witten conjecture states that $\tau$, which is a formal power series of variables
$t^0, t^1, t^2, \dots$, is uniquely determined by the following conditions:
\begin{itemize}
\item The string equation:
\[L_{-1}\tau=\sum_{p\ge0}t^{p+1}\frac{\p \tau}{\p t^p}+\frac{(t^0)^2}2\tau-\frac{\p \tau}{\p t^0}=0;\]
\item The KdV hierarchy:
\[\frac{\p U}{\p t^p}=\p_xR_p,\quad  p=0, 1, 2, \dots,\]
where $x=t^0$, $U=\p_x^2 \log\tau$, and $R_p$ are polynomials of $U, U_x, \dots, U_{(2p)x}$ which are
determined by
\begin{align*}
&R_p(0)=0,\quad  R_0=U, \\
&(2p+1)\p_x R_p=\left(2 U \p_x +U_x+\frac14 \p_x^3\right) R_{p-1}.
\end{align*}
Here we use the following notations:
\[\p_x=\frac{\p}{\p x},\quad  U_x=\p_x U,\quad \dots, \quad U_{k x}=\p_x^k U, \quad\dots.\]
\end{itemize}
The second version of the Witten conjecture states that the partition function
$\tau$ is uniquely determined by a series
of linear differential constraints:
\begin{equation}
L_m\tau=0,\quad m\ge-1, \label{vira-1}
\end{equation}
where $L_{-1}$ is given above in the string equation, and $L_m\ (m\ge0)$ are given by
\begin{align}
L_m=&\frac12\sum_{p+q=m-1}\frac{(2p+1)!!(2q+1)!!}{2^{m+1}}\frac{\p^2}{\p t^p \p t^q}\notag\\
&+\sum_{p\ge0}\frac{(2p+2m+1)!!}{2^{m+1}(2p-1)!!}(t^p-\delta_{p,1})\frac{\p}{\p t^{p+m}}
+\frac1{16}\delta_{m,0}.
\end{align}
Here $\delta_{i,j}$ is the Kronecker symbol
\[\delta_{i,j}=\left\{\begin{array}{cc} 1, & i=j, \\ 0, & i \ne j. \end{array} \right.\]
Since the operators $L_m$ satisfy the Virasoro commutation relations
\[[L_m, L_n]=(m-n) L_{m+n},\quad m, n\ge -1,\]
the constraints \eqref{vira-1} are called the Virasoro constraints of  the two-dimensional gravity.

The Witten conjecture is proved by Kontsevich in \cite{Ko}, it has inspired active researches on the following subjects:
\begin{enumerate}
\item The two-dimensional gravity can be regarded as a string theory in the zero-dimensional space-time, i.e. a point.
By considering a similar theory in space-times with richer geometric structures such as Calabi-Yau three-folds, one
obtains the Gromov-Witten (GW) invariants theory (see for example \cite{KM, RT, BF}).
When the space-time is the complex projective line $\mathbb{CP}^1$, the analogue of the Witten conjecture holds true, i.e.
the generating function of GW invariants is given by a tau function of the extended Toda integrable hierarchy, and it also satisfies the
Virasoro constraints \cite{EY, Ge, OP, DZ2, M}. For more general space-times having certain
nice properties (e.g. toric Fano), there are also some general results \cite{D1, DZ1, Gi, Te}.

\item  The two-dimensional gravity can be regarded as a field theory with spin two on Riemann surfaces.
One can also consider fields with spin $n$, then the resulting theory is just the Landau-Ginzburg (LG) model
for the $A_{n-1}$ singularity \cite{W2, PV, CH, FSZ}. By considering the LG model for general
singularities, one obtain the Fan-Jarvis-Ruan-Witten (FJRW) invariants theory (\cite{FJR1, FJR2}).
When the singularities are of ADE type, analogues of the Witten conjecture also hold true, i.e.
the generating functions of FJRW invariants are given by particular tau functions of the associated Drinfeld-Sokolov integrable hierarchies, and they also satisfy the Virasoro constraints \cite{FJR1, FJR2, FJMR, Wu, FGM}.
\end{enumerate}

In general, the first version of the Witten conjecture has the following analogue for the GW and FJRW invariants:
The generating functions of the GW and FJRW invariants should be given by particular tau functions of certain hypothetical integrable hierarchies of KdV type. Assuming the semisimplicity of the underlining Frobenius manifolds and the validity of the Virasoro constraints for the partition functions, the hypothetical integrable hierarchies can be constructed in terms of the Frobenius manifold structures, as shown in \cite{D1, DZ1}. An important question is whether there
is an analogue of the second version of the Witten conjecture for the GW and FJRW invariants, that is, whether there exist sufficiently many linear differential constraints which uniquely determine the generating functions of the GW and FJRW invariants.

The main result of the present paper is an affirmative answer to the question above for FJRW invariants of ADE singularities based on the result of Bakalov and Milanov \cite{BM}.
They proved that the partition function (also called the total descendant potential) associated to a simple singularity of type $X_\ell$ satisfies the $\W$-constraints constructed from
the $\W$-algebra of the affine Lie  algebra of type
$X_\ell^{(1)}$ and a twisted module of the corresponding vertex algebra.
Note that $\W$-algebras were first introduced by Zamolodchikov in \cite{Z} and by Fateev,
Lukyanov in \cite{FL} for the $A_n$ cases in the setting of conformal field theory.
Then Feigin and E. Frenkel \cite{FF1, FF} defined $\W$-algebras for any affine Lie
algebra as the intersection of kernels of certain screening operators associated to
the corresponding vertex algebra. 

In this paper, we prove that the $\W$-constraints constructed by Bakalov and Milanov
for simple singularities uniquely determine the partition function up to a constant factor,
as it is conjectured by Bakalov and Milanov in their paper. The answer for other cases is still unknown, and we will discuss this problem in subsequent publications.
Let us take the case of $A_1$ singularity for example to illuminate our method
to prove the uniqueness of solution of the corresponding $\W$-constraints.
In this case the $\W$-constraints are just the Virasoro constraints \eqref{vira-1}.
Note that the Virasoro constraints are linear PDEs,
so when we say that their solution $\tau$ is unique, we mean that $\tau$ is uniquely determined up to a constant factor. To prove
this assertion, we only need to show, if $\tau$, as a formal power series of $t^0, t^1, \dots$ satisfying $\tau(0)=0$, then $\tau$ vanishes itself.
To this end we introduce a degree on the ring of formal power series of $t^0, t^1, \dots$ as follows:
\[\deg t^p =p+\frac12,\]
If $\tau$ does not vanish then it must contain some nonzero monomials with the lowest degree. We denote one of them by
\[c\,t^{p_1}\dots t^{p_k},\]
where $c\ne 0$. The degree of this monomial is
\[D=\left(p_1+\frac12\right)+\cdots+\left(p_k+\frac12\right).\]
Now we take $m=p_k-1$ and denote the Virasoro constraint $L_m\tau=0$ as
\begin{align}
c_m \frac{\p \tau}{\p t_{m+1}}=&\sum_{p=0}^{m-1}a_{p,m}\frac{\p^2 \tau}{\p t^p\p t^{m-1-p}}
+\sum_{p\ge0}b_{p,m}t^p\frac{\p \tau}{\p t^{p+m}}\nn\\
&+\frac1{16}\delta_{m,0}\tau+\frac{(t^0)^2}2\delta_{m,-1}\tau \label{vira},
\end{align}
where $c_m=\frac{(2m+3)!!}{2^{m+1}}$, and $a_{p,m}$, $b_{p,m}$ are some constants.
Then the monomials with the lowest degree in the left hand side of \eqref{vira} must contain
\[c_m\, \frac{\p}{\p t^{p_k}}\left(c\, t^{p_1}\dots t^{p_k}\right)\]
which has degree
\[d_1=D-p_k-\frac12.\]
On the other hand, if we denote by $d_2$ the degree of the lowest degree monomials of the right hand side of \eqref{vira},
then we have
\[d_2\ge D-p_k+1>d_1.\]
So if  the Virasoro constraint \eqref{vira} holds true, then the coefficients of the monomials with degree $d_1$ in the left hand side must be zero.
Note that $c_m$ is always nonzero, so we must have $c=0$, and we arrive at a contradiction.

In this proof, the dilaton shift $t^1\mapsto t^1-1$ and the fact that $c_m\ne 0$ play crucial roles. In the general cases,
we also have the dilaton shift, then the key step of our proof is to show that the coefficient of a term in the $\W$-constraint
with the lowest degree does not vanish. To this end, one must check carefully the structure of these $\W$-constraints, and find out this term.

Our paper is organized as follows. In Section 2, we recall Bakalov and Milanov's construction of $\W$-constraints for simple singularities.
In Section 3, we prove our main theorem. In the last section, we  take the $D_4$ singularity as an example to show an application of our uniqueness theorem of the $\W$-constraints.

\section{Lattice vertex algebras and their twisted modules}\label{sec-2}

In this section, we recall Bakalov and Milanov's construction of $\W$-constraints of the total descendant potential associated to
a simple singularity \cite{BM}.

Let $(Q, (\ \midd\ ))$ be the root lattice of a simple Lie algebra $\g$ of ADE type. By definition, $Q$ is a free abelian group with generators $\alpha_1, \dots, \alpha_\ell$, and $(\ \midd\ )$ is a symmetric positive definite quadratic form over $Q$ such that
$a_{ij}=(\alpha_i\midd\alpha_j)$ give the entries of the Cartan matrix of $\g$.
We can define a vertex algebra associated to this lattice in the following way (see \cite{Kac} for more
details):
Let us first choose a bimultiplicative function $\e:Q \times Q \to \{\pm 1\}$ such that
\[\e(\alpha, \alpha)=(-1)^{|\alpha|^2(|\alpha|^2+1)/2},\]
where $|\alpha|^2=(\alpha\midd\alpha)$. Note that there always exists such a function. For example, we can take
\[\e(\alpha, \beta)=(-1)^{((1-\sigma)^{-1}\alpha\midd\beta)},\]
where $\sigma$ is a Coxeter transformation of $Q$. By using the bimultiplicativity property
\[\e(\alpha+\beta, \gamma)=\e(\alpha,\gamma)\e(\beta, \gamma),\quad
\e(\alpha, \beta+\gamma)=\e(\alpha,\beta)\e(\alpha, \gamma),\]
it is easy to see that $\e(\cdot, \cdot)$ satisfies the 2-cocycle condition (with trivial $Q$-action on $\{\pm 1\}$)
\[\e(\alpha,\beta)\e(\alpha+\beta,\gamma)=\e(\alpha,\beta+\gamma)\e(\beta,\gamma),\]
so we can introduce the \textit{twisted group algebra}
\[\C_\e[Q]=\mathrm{Span}_\C\{e^\alpha|\alpha \in Q\}\]
whose associative multiplication is defined by
\[e^\alpha e^\beta=\e(\alpha,\beta)e^{\alpha+\beta}.\]

Let us denote by $\h=\C\otimes_\Z Q$ the complexification of $Q$, and extend the quadratic form $(\ \midd\ )$ to $\h$ linearly.
Define the current algebra $\hh$ associated to $(\h, (\ \midd\ ))$ by
\[\hh=\h[t, t^{-1}]\oplus \C K.\]
On $\hh$ we have the following Lie algebra structure:
\[[\phi\, t^m, \phi'\, t^n]=(\phi \midd \phi')\,m\,\delta_{m+n,0}\,K,\quad  [\phi\, t^m, K]=0,\]
where $\phi,\phi'\in\h$ and $m, n\in\Z$. Introduce the \textit{bosonic Fock space}
\[\F=S(\h[t^{-1}]t^{-1}).\]
Then it is well known that $\F$ admits a level $1$ irreducible $\hh$-module structure which can be fixed by the following conditions:
\[K \cdot 1=1,\quad \phi\, t^m \cdot 1=0\ (m \ge 0),\quad \phi\, t^m \cdot s=\phi\, t^m\,s\ (m<0,\ s\in \F).\]

Now the lattice vertex algebra associated to $Q$ is defined as
\[V_Q=\F\otimes \C_\e[Q]\]
with vacuum vector $1 \otimes e^0$. To define the state-field correspondence $Y$, we need to introduce some endomorphisms on $V_Q$ as follows:
\begin{itemize}
\item For $\phi\, t^m \in \hh$, we define
\[\phi\, t^m(s\otimes e^\alpha)=\left(\phi\, t^m \cdot s\right) \otimes e^\alpha+\delta_{m,0}\,(\phi\midd\alpha)\,s\otimes e^\alpha.\]

\item For $e^\beta\in \C_\e[Q]$, we define
\[e^\beta(s\otimes e^\alpha)=\e(\beta, \alpha)s\otimes e^{\beta+\alpha}.\]

\item For $\gamma\in Q$ and an indeterminant $z$, we define
\[z^\gamma(s\otimes e^\alpha)=z^{(\gamma\midd\alpha)}s\otimes e^\alpha.\]
\end{itemize}
Then the state-field correspondence $Y$ of $V_Q$ can be generated by the following vertex operators:
\begin{align*}
&\phi(z)=Y(\phi\, t^{-1}\otimes e^0,z)=\sum_{n\in \Z} (\phi\, t^n)z^{-n-1},\\
&Y_\alpha(z)=Y(1\otimes e^\alpha, z)=e^\alpha\, z^\alpha\,
\exp \left(\sum_{n\ge1}(\alpha\, t^{-n})\frac{z^n}{n}\right)
\exp \left(\sum_{n\ge1}(\alpha\, t^{n})\frac{z^{-n}}{-n}\right).
\end{align*}

It's easy to see that $\F$ is a vertex subalgebra of $V_Q$.
\begin{dfn}
The $\W$-algebra associated to $Q$ is the vertex subalgebra of $\F$ defined by
\[\W=\{s\in\F\mid {e^{\alpha_i}}_{(0)}(s)=0,\ i=1, \dots, \ell\},\]
where ${e^\alpha}_{(0)}=\mathrm{Res}_{z=0} Y_\alpha(z)$ is called the screening operator of $\alpha \in Q$.
\end{dfn}

Let us introduce a degree on $\F$:
\begin{equation}
\deg \phi\, t^{-n}=n,\ \forall\ \phi\in\h,\ n=1, 2, \dots, \label{deg-1}
\end{equation}
then B. Feigin, E. Frenkel and I. Frenkel proved the following theorem.

\begin{thm}[\cite{FF, FI}]\label{thm-2-2}
The $\W$-algebra is generated by some elements $w_1, \dots, w_\ell$ of $\F$ such that $\deg w_i=d_i=m_i+1\ (i=1, \dots, \ell)$,
where $1=m_1\le m_2 \le \dots \le m_\ell$ are exponents of the Weyl group $W$ of $Q$. In particular, let
$I_i(\alpha_1, \dots, \alpha_n)\ (i=1, \dots, \ell)$ be the generators of $S(\h)^W$, then we can choose $w_i$ as
\begin{equation}
w_i=I_i(\alpha_1 t^{-1}, \dots, \alpha_\ell t^{-1})+J_i, \label{WIJ}
\end{equation}
where $J_i$ belong to the ideal in $\F$ that is generated by $\phi\, t^{-n}$ for all $\phi\in\h$ and $n\ge2$.
\end{thm}

In order to obtain the $\W$-constraints from the $\W$-algebra, we need to consider a twisted module $M$ of the vertex algebra $\F$.
Let $\sigma$ be a Coxeter transformation on $\h$. We denote by $h$ the order of $\sigma$, which is called the Coxeter number of $W$.
Then $\sigma$ has eigenvalues $\zeta^{m_j}$, where $\zeta=\exp(2\,\pi\,\sqrt{-1}/h)$, $j=1, \dots, \ell$.
Suppose $\phi^i$ and $\phi^j$ are eigenvectors of $\sigma$ with eigenvalues $\zeta^{m_i}$ and $\zeta^{m_j}$ respectively,
then we have
\[(\phi^i\midd\phi^j)=(\sigma(\phi^i)\midd\sigma(\phi^j))=\zeta^{m_i+m_j}(\phi_i\midd\phi_j),\]
which implies that $(\phi^i\midd\phi^j)\ne 0$ if and only if $i+j=\ell+1$. So we can choose a basis $\{\phi^i\}$ of $\h$
consisting of eigenvectors of $\sigma$ such that
$(\phi^i\midd\phi^j)=\delta_{i+j,\ell+1}$.
We also use the notation $\phi_j=\phi^{\ell+1-j}$ to denote the dual basis $\{\phi_i\}$ of $\{\phi^i\}$, and we denote $\eta_{ij}=\delta_{i+j,\ell+1}$.

The twisted module $M$ is defined as the following polynomial ring
\[M=\C[q^{i,p}\midd i=1, \dots, \ell;\ p=0, 1, 2, \dots],\]
and the action of $\F$ is given by the map
\[a\in\F \mapsto Y^M(a,\lambda)=\sum_{n\in\frac1h\Z}a_{(n)}\lambda^{-n-1},\quad \mbox{where } a_{(n)}\in\mathrm{End}(M).\]
Here the twisted state-field correspondence $Y^M$ can be generated by the following simple operators and certain reconstruction theorems of twisted modules of vertex algebras (see \cite{BK, BM} for more details):
\begin{itemize}
\item If $a=1$, then $Y^M(1, \lambda)=\mathrm{id}_M$.
\item If $a=\phi^jt^{-1}$, then
\[Y^M(\phi^jt^{-1}, \lambda)=\sum_{p\in\Z} \phi^j_{(p+\frac{m_j}{h})}\lambda^{-1-p-\frac{m_j}{h}},\]
where
\begin{equation}
\phi^j_{(p+\frac{m_j}{h})}=\frac{\Gamma(\frac{m_j}{h}+p+1)}{\Gamma(\frac{m_j}{h})}\frac{\p}{\p q^{j,p}}, \quad
{\phi_j}_{(-p-\frac{m_j}{h})}=\frac{\Gamma(\frac{m_j}{h})}{\Gamma(\frac{m_j}{h}+p)}q^{j,p} \label{boson-1}
\end{equation}
for $p=0, 1, 2, \dots$, and $j=1, \dots, \ell$.
\end{itemize}

To write down the vertex operators for other types of elements of $\F$
we need to introduce some notations.
First, for a twisted field
\[A(\lambda)=\sum_{n\in\frac1h\Z}a_{(n)}\lambda^{-n-1}\]
we denote $A(\lambda)=A(\lambda)_++A(\lambda)_-$, where
\[
A(\lambda)_+=\sum_{n\in\frac1h\Z,\ n\ge0}a_{(n)}\lambda^{-n-1},\quad
A(\lambda)_-=\sum_{n\in\frac1h\Z,\ n < 0}a_{(n)}\lambda^{-n-1}.
\]
Then for two twisted fields $A(\mu), B(\lambda)$ we define
\[:A(\mu)\,B(\lambda):=A(\mu)_-\,B(\lambda)+B(\lambda)\,A(\mu)_+, \quad \langle A(\mu)\,B(\lambda)\rangle=[A(\mu)_+, B(\lambda)_-].\]
For more fields, e.g. $A, B, C$, we assume that
\[:A\,B\,C:=:A\left(:B\,C:\right): \]
and so on.

Next, let us define
\[P^{ij}(\mu,\lambda)=\langle Y^M(\phi^it^{-1},\mu)Y^M(\phi^jt^{-1},\lambda)\rangle.\]
One can show that
\[P^{ij}(\mu,\lambda)=\eta_{ij}\,\mu^{-\frac{m_i}{h}}\lambda^{-\frac{m_j}{h}}
\frac{\frac{m_i}{h}\mu+\frac{m_j}{h}\lambda}{(\mu-\lambda)^2},\quad \mbox{for } |\lambda|<|\mu|.\]
Denote $s=\mu-\lambda$, then it is easy to see that
\[P^{ij}(\mu,\lambda)=\frac{\eta_{ij}}{s^2}+\sum_{k=0}^\infty P^{ij}_k(\lambda)s^k,\]
where
\[P^{ij}_k(\lambda)=\eta_{ij}(-1)^k\left(1-\frac{m_i}{h}\right)\frac{\Gamma\left(\frac{m_i}{h}+k+1\right)}{k!(k+2)\Gamma\left(\frac{m_i}{h}\right)}
\lambda^{-k-2}.\]

Now we can define $Y^M(a, \lambda)$ for a general $a\in\F$. Suppose $a\in\F$ is a monomial of the form
\[a=\phi^{\alpha_1}t^{-k_1-1}\otimes\cdots\otimes\phi^{\alpha_r}t^{-k_r-1},\]
then the vertex operator $Y^M(a,\lambda)$ can be obtained from the Wick theorem as follows:
\begin{equation}
Y^M(a,\lambda)=\sum_J\left(\prod_{(i,j)\in J}\p_\lambda^{(k_j)}P^{ij}_{k_i}(\lambda)\right)
:\left(\prod_{l\in J'}\p_\lambda^{(k_l)}Y^M(\phi^{\alpha_l}t^{-1},\lambda)\right):. \label{YM}
\end{equation}
Here the summation is taken over all the collections $J$ of disjoint ordered pair $(i_1, j_1), \dots, (i_s, j_s)\subset \{1, \dots, r\}$
such that $i_1<\cdots<i_s$ and $i_l<j_l$ for any $l=1, \dots, s$, and $J'=\{1, \dots, r\}\backslash J$, $\p_\lambda^{(k)}=\p_\lambda^k/k!$.
Note that the set $J$ or $J'$ can be an empty set, in such case the corresponding product is set to be $1$.

Suppose $w$ is an element of $\W$, then by using the fact that $w$ is $W$ invariant one can show that
\[Y^M(w, \lambda)=\sum_{m\in\Z}w_{(m)}\lambda^{-m-1}.\]
Let $w_1, \dots, w_\ell$ be a set of generators of $\W$, we denote
\[W_{i,m}=\mathrm{Res}_{\lambda=0} \left(\lambda^m Y^M(w_i,\lambda)\right),\quad \mbox{where } m\in\Z.\]
These operators are called the $\W$ operators associated to $Q$.

We introduce the dilaton shift $t^{i,p}=q^{i,p}+\delta_{i,1}\delta_{p,1}$, and  complete $M$ to
\[\hat{M}=\C[[t^{i,p}\midd i=1, \dots, \ell;\ p=0, 1, 2, \dots]].\]
It is easy to see that $\hat{M}$ is also a twisted module of $\F$.

The main result of \cite{BM} can be stated as follows: 
\begin{thm}[\cite{BM}]
Let $\tau$ be the total descendant potential of the semisimple Frobenius manifold associated to a simple singularity of type $X_\ell$,
and $W_{i,m}\ (i=1, \dots, \ell;\ m\in\Z)$ be the $\W$ operators associated to the root lattice of type $X_\ell$. Then as an element of $\hat{M}$ the function $\tau$ satisfies the $\W$-constraints
\begin{equation}
W_{i,m}\tau=0,\quad i=1, \dots, \ell,\ m\ge0. \label{w-constraint}
\end{equation}
\end{thm}

In this paper, we are to prove the  following theorem.
\begin{thm}[Main Theorem]\label{main-thm}
The solution to the $\W$-constraints given in \eqref{w-constraint} is unique up to a constant factor.
\end{thm}

\section{Proof of the Main Theorem}

In this section we proof the main theorem (Theorem \ref{main-thm}) of this paper.
Our method is similar to the one given in Section 1 for the $A_1$ singularity.  
We first assume $\tau(0)=0$ and $\tau \not\equiv 0$. Introduce a gradation on $\hat{M}$ by defining
\[\deg t^{i,p}=p+\frac{m_i}{h},\]
and consider the nonzero monomials of $\tau$ with the lowest degree. There may exist several monomials having the lowest
degree. We choose an arbitrary one, say,
\begin{equation}
c\,t^{i_1, p_1}\cdots t^{i_k,p_k}. \label{term}
\end{equation}
Then we will consider the equation $W_{i_k,p_k} \tau =0$.
By separating the nonzero monomials of the left hand side of this equation with the lowest degree, we can show that $c=0$ and thus arrive at a contradiction
with the assumption $\tau\not\equiv 0$.

\begin{lem}
If we assume
\begin{equation}
\deg q^{i,p}=p+\frac{m_i}{h}, \label{deg}
\end{equation}
then $W_{i,m}$, as an endomorphism of $M$, has degree $m_i-m$, i.e. $W_{i,m}(M^d)\subset M^{d+m_i-m}$, where $M^d$ is the homogeneous
component of $M$ of degree $d$ with respect to the gradation \eqref{deg}.
\end{lem}
\begin{prf}
By definition, it is easy to see that
\[\deg \phi^j_{(p+\frac{m_j}{h})}=-p-\frac{m_j}{h}, \quad \mbox{for } p\in \Z,\]
so if we assume $\deg \lambda=-1$, then we have
\[\deg Y^M(\phi\,t^{-1}, \lambda)=1,\quad \mbox{for } \phi\in\h,\]
which coincides with the degree of $\phi\,t^{-1}$ (see \eqref{deg-1}).

Note that $\deg P^{ij}_k(\lambda)=k+2$, so we have
\[\deg \p_{\lambda}^{(k_j)}P^{ij}_{k_i}(\lambda)=(k_i+1)+(k_j+1), \quad
\deg \p_{\lambda}^{(k_l)}Y^M(h t^{-1},\lambda)=k_l+1.\]
Thus it follows from the definition \eqref{YM} that 
\[\deg Y^M(a, \lambda)=\deg a, \quad \mbox{for all homogeneous } a\in \F.\]
In particular,
\[\deg W_{i,m}=\deg Y^M(w_i, \lambda)-m-1=m_i-m.\]
The lemma is proved.
\end{prf}

The above lemma shows that if we do not perform the dilaton shift, then $W_{i,m}$ are homogeneous. On the other hand, the dilaton shift
gives terms with lower degrees. Note that our aim is to find out the terms in $W_{i,m}$ which have the lowest degree, so we only need
to consider the terms that contain $t^{1,1}$.

\begin{lem}\label{lem-2}
The operator $W_{i,m}\ (m\ge0)$ has the following expression:
\begin{equation}
W_{i,m}=\sum_{d=0}^{m_i}\left(t^{1,1}-1\right)^d W_{i,m}^{(d)}, \label{form}
\end{equation}
where $W_{i,m}^{(d)}$ are differential operators on $\hat{M}$ whose coefficients do not depend on $t^{1,1}$.
\end{lem}
\begin{prf}
Every monomial of $W_{i,m}$ can be written as composition of bosons \eqref{boson-1}, and the number of composed bosons is less
or equal to the degree of $w_i$, i.e. $d_i=m_i+1$. So we have the following expression
\[W_{i,m}=\sum_{d=0}^{d_i}\left(t^{1,1}-1\right)^d W_{i,m}^{(d)},\]
where $W_{i,m}^{(d_i)}$ is a constant. If we replace $(t^{1,1}-1)$ by $q^{1,1}$ and consider the gradation \eqref{deg}, then
\[\deg \left((t^{1,1}-1)^{d_i}W_{i,m}^{(d_i)}\right)=\left(1+\frac1h\right)(m_i+1)>m_i-m=\deg W_{i,m},\]
which implies $W_{i,m}^{(d_i)}=0$. The lemma is proved.
\end{prf}

\begin{lem}\label{lem-3}
With an appropriate choice of  the generators $w_1, \dots, w_\ell$ of the $\W$-algebra
we have
\begin{equation}
W_{i,m}^{(m_i)}=c_{i,m}\frac{\p}{\p t^{i,m}}, \label{leading}
\end{equation}
where $c_{i,m}$ are nonzero constants.
\end{lem}

In order to prove the above lemma, we first need to prove the following lemma.

\begin{lem}\label{lem-4}
The generators $I_1,\dots, I_\ell$ of $S(\h)^W$ can be chosen to have the 
form
\begin{equation}
I_i(\phi_1, \dots, \phi_\ell)=\phi_1^{m_i}\phi_{\ell+1-i}+\sum_{d=0}^{m_i-1}\phi_1^d I_i^{(d)}(\phi_2, \dots, \phi_\ell),
\quad i=1, \dots, \ell.
\end{equation}
Here the dual basis $\{\phi_1, \dots, \phi_\ell\}$ are regarded as coordinates of $\h$ with respect to
the basis $\{\phi^1, \dots, \phi^\ell\}$ (see Section \ref{sec-2}), 
and the invariant polynomials $I_1, \dots, I_\ell$ are represented as polynomials 
of $\phi_1, \dots, \phi_\ell$.
\end{lem}
\begin{prf}
We first take an arbitrary set $\{I_1, \dots, I_\ell\}$ of generators of $S(\h)^W$.
Let $d_i' (1\le i\le \ell) $ be the degree of $I_i$ which is not necessarily equal to $m_i+1$.

Let $J$ be the Jacobian determinant of $I_1, \dots, I_\ell$ with respect to $\phi_1, \dots, \phi_\ell$. Then $J$ is a constant multiple of
the product of linear equations of all the walls of a Weyl chamber. Note that $\phi^1=(1, 0, \dots, 0)$, which is an eigenvector of the
Coxeter transformation $\sigma$, does not lie on any wall of a Weyl chamber, so we have $J(1, 0, \dots, 0)\ne 0$.

The determinant $J$ contains $\ell!$ summands, in which there is at least one that does not vanish at $(1, 0, \dots, 0)$.
We can re-number $I_1, \dots, I_\ell$ such that this summand is given by
\[(-1)^{\frac{\ell(\ell-1)}{2}}\frac{\p I_1}{\p \phi_\ell}\frac{\p I_2}{\p \phi_{\ell-1}}\dots\frac{\p I_\ell}{\p \phi_1},\]
then we have $\frac{\p I_i}{\p \phi_{\ell+1-i}}(1, 0, \dots, 0)\ne 0$ for $i=1, \dots, \ell$.

By using the above fact, we can rescale $I_1, \dots, I_\ell$ such that
\[I_i=\phi_1^{d_i'-1}\phi_{\ell+1-i}+\mbox{other terms}.\]
From the fact that $\sigma(\phi_i)=\zeta^{-m_i}\phi_i$, $m_1=1$ and $\sigma(I_i)=I_i$, it follows that $d_i'$ must be equal to $m_i+1$, and the ``other terms" in the above expression of $I_i$ only contains monomials with degrees less than $d_i'-1$
(in the $D_{2n}$ case, a linear recombination of $I_{n}$ and $I_{n+1}$ may be needed). The lemma is proved.
\end{prf}

\begin{rmk}
The proof of the above lemma follows the one given for Theorem 3.19 in \cite{Hum}.
\end{rmk}

\begin{prfn}{Lemma \ref{lem-3}}
By counting degrees and the number of bosons, it is easy to see that $W_{i,m}^{(m_i)}$ must take the form
\eqref{leading} (in the $D_{2n}$ case, a linear recombination of $I_{n}$ and $I_{n+1}$ may be needed).
We only need to prove that $c_{i,m}\ne 0$.

Let us call $(t^{1,1}-1)^{m_i}W_{i,m}^{(m_i)}$ the leading term of $W_{i,m}$. This leading term is the composition of $d_i$ bosons.
Note that if we split $w_i$ into $I_i+J_i$ (see \eqref{WIJ}), then the leading term must come from $I_i$ since $J_i$ contains less
bosons. Similarly, suppose $a$ is a monomial of $I_i$, we consider the operator $Y^M(a, \lambda)$ (see \eqref{YM}), then the leading
term must come from the summands with $J=\emptyset$ since other summands contain less bosons. So in order to find out the leading term
we only need to investigate the invariant polynomial $I_i$, and we can omit all terms that contain $P^{ij}_k(\lambda)$ in \eqref{YM}.

According to Lemma \ref{lem-4}, the leading term of $W_{i,m}$ is the same with the one of
\[\mathrm{Res}_{\lambda=0}\left(\lambda^m Y^M((\phi_1\,t^{-1})^{m_i}(\phi^i\,t^{-1}), \lambda)\right).\]
By definition,
\begin{align*}
Y^M(\phi_1\,t^{-1}, \lambda)&=h\,(t^{1,1}-1)\lambda^{\frac1{h}}+\dots,\\
Y^M(\phi^i\,t^{-1}, \lambda)&=\frac{\Gamma(\frac{m_i}{h}+m+1)}{\Gamma(\frac{m_i}{h})}\frac{\p}{\p t^{i,m}}\lambda^{-1-m-\frac{m_i}{h}}+\dots,
\end{align*}
so the leading coefficient $W_{i,m}^{(m_i)}$ reads
\[W_{i,m}^{(m_i)}=h^{m_i}\frac{\Gamma(\frac{m_i}{h}+m+1)}{\Gamma(\frac{m_i}{h})}\frac{\p}{\p t^{i,m}}.\]
The lemma is proved.
\end{prfn}

Now we are ready to prove Theorem \ref{main-thm}. From Lemma \ref{lem-2} and Lemma \ref{lem-3} it follows that the lowest degree term in $W_{i,m}$ is given by
\[(-1)^{m_i}c_{i,m}\frac{\p}{\p t^{i,m}}.\]
It's action on the term \eqref{term} implies
\[(-1)^{m_{i_k}}c_{i_k,p_k}\frac{\p}{\p t^{i_k,p_k}}\left(c\,t^{i_1, p_1}\cdots t^{i_k,p_k}\right)=0\]
so $c=0$. The theorem is proved.

\section{Application of $\W$-constraints}

If we begin with $\tau(0)=1$ instead of $\tau(0)=0$, then the proof given in the last section provides an algorithm to compute the partition function $\tau(t)$.
In particular, when the singularity $X_\ell$ is of $A$ type, the explicit form of $\W$-constraints was obtained in \cite{FL, Go, AvM}, and the corresponding partition functions $\tau$ can be obtained by solving these $\W$-constraints recursively. In this section, we consider the $D_4$ singularity.

We take $\h=\C^4$ with the natural basis $e_1, \dots, e_4$ first. The metric $(\cdot\midd\cdot)$ is taken as the canonical one
$(e_i\midd e_j)=\delta_{i,j}$, and the generators $\alpha_1, \dots, \alpha_4$ of the root lattice $Q$ are chosen as
\[\alpha_1=e_1-e_2,\quad \alpha_2=e_2-e_3,\quad \alpha_3=e_3-e_4, \quad \alpha_4=e_3+e_4.\]
Then the Cartan matrix $A=((\alpha_i\midd \alpha_j))$ reads
\[A=\left(\begin{array}{cccc} 2  & -1 &  0 &  0 \\ -1 &  2 & -1 & -1 \\
0  & -1 &  2 &  0 \\ 0  & -1 &  0 &  2 \end{array}\right),\]
and $Q=\Z\alpha_1\oplus\Z\alpha_2\oplus\Z\alpha_3\oplus\Z\alpha_4$ is a root lattice of $D_4$ type.

The Coxeter transformation is taken as $\sigma=R_1R_2R_3R_4$, where $R_i\ (i=1, \dots, 4)$ is the reflection with respect to $\alpha_i$.
The order of $\sigma$ is $h=6$, and the exponents read
\[m_1=1,\quad m_2=3, \quad m_3=3, \quad m_4=5.\]
The eigenvectors $\phi^1, \dots, \phi^4$ are chosen as
\begin{align*}
\phi^1=\overline{\phi^4}=&\frac1{\sqrt3}\left(\zeta^2\,e_1+\zeta\,e_2+e_3\right),\\
\phi^2=\overline{\phi^3}=&\frac1{\sqrt6}\left(-e_1+e_2-e_3+\sqrt{-3}\,e_4\right),
\end{align*}
where $\zeta=e^{\pi\sqrt{-1}/3}$, and $\overline{{\phi^3}}, \overline{{\phi^4}}$ stand for the complex conjugations of $\phi^3, \phi^4$.

Bakalov and Milanov showed in \cite{BM} that the generators of the corresponding $\W$-algebra can be chosen as
\begin{align*}
\tilde{w}_i&=\sum_{j=1}^4 e^{e_j}{}_{(-m_i-1)}\,e^{-e_j}+\sum_{j=1}^4 e^{-e_j}{}_{(-m_i-1)}\,e^{e_j},\quad i=1,2,4,\\
\tilde{w}_3&=(e_1\,t^{-1})(e_2\,t^{-1})(e_3\,t^{-1})(e_4\,t^{-1}).
\end{align*}
These generators do meet the requirement of Lemma \ref{lem-3}, so we modify them to the following set of generators:
\begin{align*}
w_1=&\tilde{w}_1, \quad w_2=3\,\tilde{w}_2+\sqrt{-3}\,\tilde{w}_3,
\quad w_3=3\,\tilde{w}_2-\sqrt{-3}\,\tilde{w}_3,\quad w_4=\tilde{w}_4.
\end{align*}
They have the explicit forms 
\begin{align*}
w_1=& \phi^{\alpha}\phi_{\alpha},\\
w_k=& 2\phi^{\alpha}\phi_{\alpha}^3+\frac34\phi^{\alpha,2}\phi_{\alpha}^2+\frac18w_1^2
+\phi^{1}(\phi_k)^2\phi^{4}+\frac16(\phi^{k})^4\\
&-\frac13\phi^{k}(\phi_k)^3-\frac{\sqrt{2}}{3}\left((\phi^{1})^3+(\phi^{4})^3\right)\phi^{k} \quad (k=2,\, 3),\\
w_4=&\frac{2}{5}\phi^{\alpha,5}\phi_\alpha+\frac14\phi^{\alpha,4}\phi_\alpha^2+\frac19\phi^{\alpha,3}\phi_\alpha^3
+\frac{(\phi^1)^6-(\phi^2)^6-(\phi^3)^6+(\phi^4)^6}{3240}\\
&-\frac{(\phi^1)^3+(\phi^4)^3}{324 \sqrt{2}}\left((\phi^2+\phi^3)^3+3(\phi^2+\phi^3)\phi^1\phi^4\right)\\
&+\frac{1}{432}\left(\phi^1\phi^4(\phi^2+\phi^3)^4+6(\phi^1\phi^4)^2(\phi^2+\phi^3)^2+\frac83(\phi^1\phi^4)^3\right.\\
&\left.+\phi^2\phi^3((\phi^2)^4-2(\phi^2)^3\phi^3-2\phi^2(\phi^3)^3+(\phi^3)^4)+\frac{10}3(\phi^2\phi^3)^3\right)\\
&+\frac1{18}\left(\phi^{1,3}\phi^4+\phi^{4,3}\phi^1\right)\left((\phi^2+\phi^3)^2+2\phi^1\phi^4\right)\\
&-\frac{1}{27\sqrt{2}}\left(((\phi^1)^3+(\phi^4)^3)(\phi^{2,3}+\phi^{3,3})
+3(\phi^2+\phi^3)((\phi^1)^2\phi^{1,3}+(\phi^4)^2\phi^{4,3})\right)\\
&+\sum_{k=2}^3\frac{\phi^{k,3}}{18}\left(\frac{2(\phi^k)^3-(\phi_k)^3}{3}+\phi^2\phi^3(2\phi_k-\phi^k)
+2\phi^1\phi^4(\phi^2+\phi^3)\right),
\end{align*}
where $\phi^{i,n}=\phi^i\,t^{-n}$, $\phi_i^n=\phi_i\,t^{-n}$, $\phi^i=\phi^{i,1}$, $\phi_i=\phi_i^1$, and summation with
respect to the repeated upper and lower Greek index $\alpha$ is assumed.

One can obtain the $\W$ operators $W_{i,m}$ from these generators, whose explicit form is omit here.
By using the $\W$-constraints $W_{i,m}\tau=0\ (m\ge0)$, we can obtain all Taylor coefficients of $\tau(t)$. In particular,
if we consider the genus expansion
\[\F(t)=\log \tau(t)=\sum_{g=0}^\infty\F_g(t),\]
and restrict $\F_0(t)$ to the small phase space
$\mathrm{Spec}\,\C[t^{1,0},t^{2,0},t^{3,0},t^{4,0}]$, then we obtain the
following potential of the Frobenius manifold associated to the
$D_4$ singularity:
\begin{equation}\label{D4_F}
F(v)=\frac{1}{2} (v^1)^2 v^4 + v^1 v^2 v^3 + \frac{((v^2)^3+(v^3)^3) v^4}{18 \sqrt{2}} + \frac{v^2 v^3 (v^4)^3}{108}
+ \frac{(v^4)^7}{272160},
\end{equation}
where $v^i=t^{i,0}$.

There are several degrees of freedom when we choose the basis $\{\phi^1,\dots,\phi^4\}$.
If we choose a different basis, then the potential $F(v)$ will be transformed to a different form via a linear coordinate
transformation. For example, if we take
\[v^1\mapsto c^{-1}\,v^1,\quad v^4\mapsto c\,v^4,\quad F \mapsto c^{-1}\,F\]
with $c=\sqrt{18\sqrt{2}}$, then
\begin{equation}
F=\frac{1}{2} (v^1)^2 v^4+v^1v^2v^3 + (v^2)^3 v^4 + 6 v^2v^3 (v^4)^3 + (v^3)^3 v^4 + \frac{54}{35} (v^4)^7
\end{equation}
which coincides with the potential derived by Dubrovin in \cite{Du1}.

If we take
\[v^1=t_1,\quad v^2=-\frac{t_X-\sqrt{3}\,t_Y}{\sqrt{2}},\quad v^3=-\frac{t_X+\sqrt{3}\,t_Y}{\sqrt{2}},\quad v^4=t_{X^2},\]
and rescale $F$ to $F/6$, then
\begin{align*}
F=&\frac{1}{12} t_1 t_X^2-\frac{1}{4} t_1 t_{Y}^2+\frac{1}{12} t_1^2 t_{X^2}-\frac{1}{216} t_X^3 t_{X^2}
-\frac{1}{24} t_X t_{Y}^2 t_{X^2}\\
&\qquad+\frac{1}{1296} t_X^2 t_{X^2}^3 -\frac{1}{432} t_{Y}^2 t_{X^2}^3 + \frac{1}{1632960} t_{X^2}^7
\end{align*}
which coincide with the potential derived by Fan \textit{et al} in \cite{FJMR}.

\vskip 2em

\noindent\textbf{Acknowledgments.}

The authors thank Boris Dubrovin for his encouragement and many helpful discussions. 
This work is partially supported by the NSFC No. 11071135, No. 11171176 and No. 11222108, and by the Marie Curie IRSES project
RIMMP.

\end{document}